\documentclass[12pt]{article}
\usepackage{amssymb}

\usepackage{amsthm}
\usepackage{amsmath}
\usepackage{latexsym}
\usepackage{graphicx}
\newtheorem{theorem}{Theorem}%[section]
\newtheorem{lemma}{Lemma}%[section]

\begin{document}
\begin{center}
{\Large \bf Anomalous Diffusion with Periodical Initial Conditions
on Interval with Reflecting Edges
\\$^*$Tsitsiashvili G.Sh, $^{**}$Yashin A.E.\\
Institute for Applied Mathematics,\\ Far Eastern Branch of RAS\\
Russia, 690041, Vladivostok, Radio str. 7,\\ IAM FEB RAS\\
E-mail: $^*$ guram@iam.dvo.ru,\\$^{**}$ ziby@mail.primorye.ru}
\end{center}
\setcounter{equation}{0}\renewcommand{\theequation}{1.\arabic{equation}}
\section{Introduction}
%{\small
In \cite{Uchaikin} a mathematical model of an anomalous diffusion
in a space was suggested. This model origins in an investigation
of processes in complex systems with variable structure: glasses,
liquid crystals, biopolymers, proteins and etc.

In this model a coordinate of a particle has stable distribution
(not normal one). As a result a density of its distribution
function satisfies an analog of diffusion equation in which second
derivative by coordinate is replaced by partial derivative.

In this paper the anomalous diffusion with periodic initial
conditions on an interval with reflecting edges is considered.
Such problem is important for example in technical mechanics for
an analysis of fuel mixing in straight flow engine \cite{Besp}
too. As A.A. Borovkov suggested \cite{Borovkov} special
probability methods are developed to analyze such a diffusion.

Main results of the paper may be represented as follows. Suppose
that $T_n(a,r,k)$ is a characteristic time of the anomalous
diffusion with a parameter $a,\;0<a\leq 2$ on $k$-dimension
interval $[-r,\;r]^k$ with $n$-periodic (by each coordinate)
initial conditions then
$$T_1(a,r,1)=T_1(a,1,1)r^{a},\;\;\frac{T_1(a,r,1)T_1(2,1,1)}{T_1(a,1,1)T_1(2,r,1)}=r^{a-2},\;\;T_n(a,1,1)=T_1(a,1,1)/n^{a}.$$
So the anomalous diffusion on the interval $[-r,r]$ for $r>1$
works faster and for $r<1$ - slower than normal one. If $k>1$ then
it is possible to choose $n$-periodical (by each coordinate)
initial conditions so that for the normal diffusion
$$T_n(2,1,k)=T_1(2,1,1)/kn^{2}.$$
\setcounter{equation}{0}\renewcommand{\theequation}{2.\arabic{equation}}
\section{Model of Anomalous Diffusion on Straight Line}
Suppose that $y(t),\;\;t\geq 0$ is homogenous random process with
independent increments and an initial condition $y(0)=0.$ The
difference $y(t)-y(\tau),\;\;t>\tau\geq 0$ has symmetric stable
distribution on the straight line $(-\infty,\;\infty)$ with a
parameter $a,\;\;0<a\leq 2$ and a characteristic function
$$M\exp(iu[y(t)-y(\tau)])=\exp(-(t-\tau)|u|^a).$$
In \cite{Uchaikin} it was shown that distribution density
$p_t=p_t(u)$ of the process $y(t)$ in a moment $t>0$ is the
generalized solution of the following differential equation with
fractional derivatives
$$\left(\frac{\partial}{\partial t}-\frac{\partial^a}{\partial
y^a}\right)p_{t}(y)=\delta(y)\delta(t).$$ Here $\partial^a
p_{t}(y)/\partial y^a$ is $a-$fractional derivative of the
function $p_{t}(y),\;\;0<a\leq 2;\;\;\delta(y),\;\delta(t)$ are
delta-functions by variables $y,\;t$ accordingly.

\setcounter{equation}{0}\renewcommand{\theequation}{3.\arabic{equation}}
\section{Construction of Anomalous Diffusion Model\\ on the Interval
$[-1,\;1].$}

Each realization of random process $y(t),\;\;t\geq 0$ may be
considered as a curve $\Gamma$ on the plane $(y,t)$ and may be
represented in parametrical form:
$y=y(\tau),\;\;t=t(\tau)=\tau,\;\;\tau\geq 0.$ Suppose that the
curve $\Gamma$ is reflected from the lines $y=1,\;y=-1$. For this
aim represent the plane $(y,\;t)$ as a transparent and infinitely
thin sheet of paper with the curve $\Gamma.$
%\linebreak
%\begin{center}
%\includegraphics[width=6.2in,height=3.3in]{rne.bmp}
%Fig. 1. Transformation of the curve $\Gamma$ into the curve
%$\gamma$
%\end{center}
Bend this sheet of the paper along the lines $y=\pm 1,\;y=\pm
3,\ldots $ into transparent strip $-1\leq y\leq 1$ with fragments
of initial curve $\Gamma$. As a result $\Gamma$ is transformed
into the curve $\gamma$
%(see fig. 1)
:
$y=Y(\tau),\;t=t(\tau)=\tau.$ Analogously to \cite{Uchaikin} the
random process $Y(t),\;\;t\geq 0$ may be considered as a model of
the anomalous diffusion on the interval $[-1,\;1]$ with reflecting
edges. It is clear that if the curve $\Gamma$ coincides with some
straight line then the curve $\gamma$ is constructed in accordance
with the law of geometrical optics: falling angle equals to
reflecting angle.

%Нашей задачей является исследование скорости сходимости
%распределения процесса аномальной диффузии в момент $t>0$ к
%равномерному распределению на этом отрезке при $t\to\infty$ в
%зависимости от различных параметров модели.

%Данная задача возникает при изучении процессов в сложных системах
%с меняющейся структурой: стекла, жидкие кристаллы, полимеры,
%протеины, биополимеры и т.д. \cite{Uchaikin}, а также при
%исследовании перемешивания топлива в прямоточных двигателях
%\cite{Besp}.
%\setcounter{equation}{0}
\setcounter{equation}{0}\renewcommand{\theequation}{4.\arabic{equation}}
\section{Analytical Representation of Reflected Diffusion\\ Process and Its Properties}
Additionally to geometrical representation of the process
$Y(t),\;\;t\geq 0$ consider its analytical representation by the
functions $f:R\to[-2,\;2),\;g:[-2,\;2)\to[-1,\;1],$
%\begin{equation}\label{0}
$$
f(u)=(u+2)/mod\;4-2,
%\end{equation}
$$
$$
%\begin{equation}\label{0.1}
g(u)= \left\{ \begin{array}{l}
u,\;-1\leq u\leq 1,\\
2-u,\;1<u< 2,\\
-2-u,\;-2\leq u < -1, \end{array}\right.
$$
%\end{equation}
\begin{equation}\label{1}
Y(t)=g(f(y(t))).
\end{equation}
Here $u/mod\;A=A\{u/A\},\;A>0;$ $\{z\}$ is the fractional part of
a real number $z.$ It is known from group theory that for
$u,\;v\in R,\;\;A>0$
\begin{equation}\label{1+-+}
(u/mod\;A)/mod\;A=u/mod\;A,\;\;(u+v)/mod\;A=(u/mod\;A+v/mod\;A)/mod\;
A.
%/\mod\;A
\end{equation}
Define on the half-interval $[-2,\;2)$ binary operation $"\oplus"$
and unary operation of an inverse $"\ominus":$
$$u\oplus v=f(u+v),\;\;\ominus u=f(-u),\;\;u,\;v\in [-2,\;2).$$
Prove that these operations on the half-interval $[-2,\;2)$ form
commutative group $C$ with the summation $\oplus$ and the inverse
$\ominus$ and the unit $0:$
$$u\oplus v=v\oplus u,\;\;u\oplus (\ominus u)=0,\;\;(u\oplus v)\oplus w=u\oplus(v\oplus w),\;\;u,v,w\in C.$$
Really using (\ref{1+-+}) obtain
$$u\oplus v=f(u+v)=(u+v+2)/mod\; 4-2=(v+u+2)/mod\; 4-2=f(v+u)=v\oplus u,$$
$$u\oplus (\ominus u)=(u+(\ominus u)+2)/mod\; 4-2=(u+(-u+2)/mod\; 4-2+2)/mod\; 4-2=$$$$=(u-u+2)/mod\; 4-2=f(0)=0,$$
$$
(u\oplus v)\oplus
w)=[f(u+v)+w+2]/mod\;4-2=[(u+v+2)/mod\;4-2+w+2]/mod\;4-2=
$$$$=(u+v+w+2)/mod\;4-2=[u+(v+w+2)/mod\;4-2+2]/mod\;4-2=$$$$
=[u+f(v+w)+2]/mod\;4-2=u\oplus(v\oplus w).$$

Prove that $f:(R,+)\to(C,\oplus)$ is a homomorphism of the
additive group $R$ of real numbers onto the group $C:$
$$f(u+v)=f(u)\oplus f(v),\;\;f(-u)=\ominus f(u).$$
Using (\ref{1+-+}) obtain $$f(u)\oplus
f(v)=f(f(u)+f(v))=[(u+2)/mod\;4-2+(v+2)/mod\;4-2+2]/mod\;4-2=$$
$$=(u+v+2)/mod\;4-2=f(u+v),$$
the equality $f(-u)=\ominus f(u)$ is a corollary of the inverse
$\ominus$ definition.
%\begin{center}
%\includegraphics[width=6.8in,height=2.4in]{rn1.bmp}
%Fig. 2. Graphics of the functions $f,\;g$
%\end{center}
\setcounter{equation}{0}\renewcommand{\theequation}{5.\arabic{equation}}
\section{Reflection Formula for Diffusion Process}
Denote by $p_t=p_t(u),\;\pi_{t}=\pi_{t}(u),\;P_{t}=P_{t}(u)$
distribution densities of random variables (r.v.)
$y(t),\;Y(t),\;f(y(t))$ accordingly. Using the formula (\ref{1})
and the graphic of the function $f$
%(see fig. 2)
find:
\begin{equation}\label{7}
P_{t}(u)=\sum_{v:\;f(v)=u}p_t(v)=\sum_{k=-\infty}^{\infty}p_t(u-4k),\;\;u\in[-2,\;2)
;\;\;P_{t}(u)=0,\;\;u\not\in[-2,\;2).
\end{equation}
It is easy to prove from \cite[chapt. 17, \S6, lemma 1]{1Fel} that
the row (\ref{7}) converges.

Define auxiliary function
$$\overline{P}_{t}(u)=\sum_{k=-\infty}^{\infty}p_t(u-4k),\;\;-\infty<u<\infty.$$
The function $\overline{P}_t(u)$ coincides with $P_{t}(u)$ for
$u\in[-2,\;2)$ and has the period $4$ and is symmetric. So from
the formula (\ref{1}) and the graphic of the function $g$
%(see
%fig. 2)
obtain
\begin{equation}\label{5.1}
\pi_{t}(u)=\overline{P}_t(u)+\overline{P}_t(u+2)=\sum_{k=-\infty}^{\infty}p_t(u-2k),\;\;u\in[-1,\;1],\;\;
\pi_{t}(u)=0,\;\;u\not\in[-1,\;1].
\end{equation}
Suppose that r.v. $z$ is distributed on $[-1,\;1],$ its density
$\mu(v)$ has continuous derivative,
$${d\mu(v)\over dv}=0,\;\;v=\pm 1$$
and r.v. $z$ doesn`t depend on random process $y(t)$ then the
formula (\ref{7}) allows to calculate distribution density
$\pi_{t}^z(u)$ of the random process $Z(t)=g(f(y(t)-1+z)),\;t\geq
0$ in the moment $t>0:$
\begin{equation}\label{5.1+}
\pi_{t}^z(u)=\int_{-1}^1(\overline{P}_t(u-v)+\overline{P}_t(u+2+v))\mu(v)dv.
\end{equation}
Remark that the formulas (\ref{5.1}), (\ref{5.1+}) which give
distribution of diffusion reflected process are analogous to
reflection method formulas which give solution of wave equation
for finite string with fixed edges \cite[chapt. III, \S 13, points
5, 6]{Vladimirov}.
\setcounter{equation}{0}\renewcommand{\theequation}{6.\arabic{equation}}
\section{Convergence of Reflected Diffusion Process\\ Distribution to Uniform Distribution}
As $f$ is the homomorphism of the additive group $R$ onto the
group $C,$ then for $0<t_1<t_2$
$$f(y(t_2))=f(y(t_1))\oplus f(y(t_2)-y(t_1)).$$
So random process $f(y(t)),\;t\geq 0$ with independent and stable
distributed increments is homogenous Markov process with state set
$[-2,\;2)$ and symmetric distribution density of a transition from
the state $u$ into the state $v$ during the time $t$
%\begin{equation}\label{3.0}
$$q_{t}(u,v)=\overline{P}_{t}(v-u),\;\;-2\leq u,\;v< 2.$$
%\end{equation}
Denote $Q_t=\inf\{q_{t}(u,v),\;\;-2\leq u,\;v< 2\}$ then from
(\ref{7}) obtain
\begin{equation}\label{3}
0<Q_t<{1\over 4},\;\;t>0.
\end{equation}
Suppose that
\begin{equation}\label{3+}
P=P(u)=\left\{
\begin{array}{l}
1/4,\;u\in[-2,\;2),\\
0,\;u\not\in[-2,\;2),
\end{array}\right.
\end{equation}
$$\pi=\pi(u)=\left\{
\begin{array}{l}
1/2,\;u\in[-1,\;1],\\
0,\;u\not\in[-1,\;1].
\end{array}\right.
$$
For function $\varphi$ on $(-\infty,\;\infty)$ define the norm
$||\varphi||=\sup\{|\varphi(u)|,\;-\infty< u<\infty\}$ and denote
$[z]$ the integer part of a real number $z.$
\begin{lemma}\label{lem1}
For arbitrary $h>0,\;\;t\geq h:$
\begin{equation}\label{6+-+}
||\pi_{t}-\pi||\leq 2(1-4Q_h)^{k-1}||P_{h}-P||\;,\;\;k=[t/h].
\end{equation}
\end{lemma}

\proof As transition density $q_t(u,v)$ is symmetric then the
formula (\ref{3+}) gives
\begin{equation}\label{3.10}
\int_{-2}^2P(v)q_t(v,u)dv={1\over
4}\int_{-2}^2q_t(u,v)dv=P(u),\;\;-2\leq u< 2.
\end{equation}
If $\Delta q_t=q_t-Q_t\geq 0$ then the formulas (\ref{3}),
(\ref{3.10}) give
$$||P_{t+1}-P||=\sup\left(\left|\int_{-2}^2P_1(v)q_t(v,u)dv-\int_{-2}^2P(v)q_t(v,u)dv\right|,\;-2\leq
u< 2\right)=$$
$$=\sup\left(\left|\int_{-2}^2P_1(v)\Delta
q_t(v,u)dv-\int_{-2}^2P(v)\Delta q_t(v,u)dv\right|,\;-2\leq u<
2\right)$$ and for $t\geq 0$
\begin{equation}\label{3.10++}
||P_{t+1}-P||\leq ||P_1-P||\sup\left(\int_{-2}^2\Delta
q_t(v,u)dv,\;-2\leq u< 2\right)=||P_1-P||(1-4Q_t).
\end{equation}
The formulas (\ref{3.10}), (\ref{3.10++}) allow to obtain by
induction
\begin{equation}\label{4}
||P_{k}-P||\leq (1-4Q_1)^{k-1}||P_{1}-P||\;,\;\;k\geq 1.
\end{equation}
Using the formulas (\ref{3}), (\ref{3.10++}) it is easy to
generalize the inequality (\ref{4}) onto arbitrary $t\geq 1$
\begin{equation}\label{4+}
||P_{t}-P||\leq (1-4Q_1)^{k-1}||P_{1}-P||\;,\;\;k=[t].
\end{equation}
The formulas (\ref{5.1}), (\ref{4+}) give
\begin{equation}\label{6}
||\pi_{t}-\pi||\leq2||P_{t}-P||\leq
2(1-4Q_1)^{k-1}||P_{1}-P||\;,\;\;k=[t],\;\;t\geq 1.
\end{equation}
The inequality (\ref{6}) may be rewritten for arbitrary positive
$h$ in the form (\ref{6+-+}). \endproof
\setcounter{equation}{0}\renewcommand{\theequation}{7.\arabic{equation}}
\section{Diffusion on Interval $[-r,\;r]$}
Suppose that $m>0,$ $r=m^{-1/a}$ and consider Markov process
$Y_r(t),$
%\\
$t\geq 0:$
%\begin{equation}\label{6.0}
$$Y_r(t)=r g(f(y(t)/r)).$$
%\end{equation}
The process $Y_r(t),\;t\geq 0 $ is obtained from the process
$y(t),\;t\geq 0$ by reflections from edges of the interval
$[-r,\;r].$ Denote by $\pi_{t,\;r}=\pi_{t,\;r}(u)$ distribution
density of r.v. $Y_r(t).$ It is obvious that $\pi_{t,\;1}=\pi_t.$
Introduce normalized r.v.
$$W_{t,\;r}={Y_r(t)\over r}=g(f(y(t)/r)),\;t\geq 0$$
with distribution density $\psi_{t,\;r}(u)=r\pi_{t,\;r}(ru).$
\begin{lemma}\label{lem2}
For $K=[tm],\;\;tm\geq 1$
\begin{equation}\label{6.1}
||\pi_{t,\;r}(u)-\pi(u/ r)/r||\leq
{2(1-4Q_1)^{K-1}||P_{1}-P||\over r}.
\end{equation}
\end{lemma}

\proof The definition of stable distribution (see \cite[chapt. 17,
\S5]{1Fel}) gives that for $t>0$ r.v.`s $y(tm),\;y(t)/r$ coincide
by distribution. So r.v.`s
%$$W_{t,\;r}=g(f(y(t)/r)),$$ причем в соответствии с (\ref{6.0})
$W_{t,\;r},\;g(f(y(tm)))$ coincide by distribution too and from
(\ref{6})
$$
||\psi_{t,\;r}-\pi||=||\pi_{tm}-\pi||\leq
2(1-4Q_1)^{K-1}||P_{1}-P||\;,\;\;K=[tm],\;\;tm\geq 1.
$$
The formula (\ref{6.1}) is proved.
\endproof
The formula (\ref{6.1}) may be rewritten conditionally as
\begin{equation}\label{6.2} T_1(a,r,1)=T_1(a,1,1) r^a.
\end{equation}

%В заключение данного раздела следует отметить, что соотношение
%вида (\ref{6.2}) может быть получено для многомерного
%диффузионного процесса с отражениями от границы множества
%достаточно общего вида, используя дифференциальные или
%интегро-дифференциальные уравнения для плотности распределения
%этого процесса. Однако в одномерном случае достаточно ограничиться
%приведенным выше рассмотрением.
\setcounter{equation}{0}\renewcommand{\theequation}{8.\arabic{equation}}
\section{Numerical Experiment}
Here analytical results of last paragraphs are compared with
results of a numerical experiment. In this experiment a closeness
of probability densities $\pi_t,\;\pi$ is investigated  for
different $t>0.$ For this aim we imitate independent r.v.`s which
coincide with r.v. $Y(t)$ and so with r.v. $g(f(t^{1/a}y(1)))$ by
distribution. R.v. $y(1)$ is imitated  approximately by normalized
sum
$$\hat{y}(1)={v_1+\ldots+v_N\over N^{1/a}}$$
of independent identically distributed r.v. $v_1,\ldots,v_N,$
$$P(v_1>t)={t^{-a}\over 2},\;\;P(v_1<-t)={|t|^{-a}\over 2},\;\;t>1.$$
Using $M$ independent realizations of r.v.
$g(f(t^{1/a}\hat{y}(1)))$ it is possible to calculate frequences
$S_j(t),\;j=0,\ldots,9$ of these realizations scorings into the
sets
$$\left[-1+{2j\over 10},\;-1+{2(j+1)\over 10}\right),\;\;j=0,\ldots,8,\;\;\left[-1+{2j\over 10},\;1\right],\;\;j=9.$$
We calculate quantities
$$S(t)=10M\sum_{j=0}^9\left(S_j(t)-{1\over 10}\right)^2$$
which characterize (analogously to $\chi^2$-statistics) deviations
of $Y(t)$ distribution densities from uniform density for
different $t>0.$

\vskip 0.5cm
\begin{center}
\setlength{\tabcolsep}{9pt}
\begin{tabular}{|c|c|c|c|}\hline
$t$  &  $S$($t$) for $a$=1.9 &
$S$($t$) for $a$=1.95  &  $S$($t$) for $a$=1.99 \\
\hline
0,01 & 9169.31 & 7867.49 & 6584.26  \\
\hline
0,02 & 4055.39 & 3053.09 & 2537.75 \\
\hline
0,03 & 2200.99 & 1298.15 & 1017.17 \\
\hline
0,04 & 1061.13 & 609.86 & 347.35 \\
\hline
0,05 & 488.14 & 289.15 & 185.43 \\
\hline
0,06 & 220.69 & 177.32 & 122.32\\
\hline
0,07 & 135.03 & 52.97 & 39.26 \\
\hline
0,08 & 102.79 & 39.37 & 13.26 \\
\hline
0,09 & 54.89 & 14.55 & 12.78 \\
\hline
\end{tabular}
\end{center}
\begin{center}
Table 1
\end{center}
Results of these calculations represented in the table 1 show how
fast do distributions of $Y(t)$ converge to uniform if the time
$t$ increases. Qualitative coincidence of numerical experiment
results with the formula (\ref{6+-+}) is demonstrated.
\setcounter{equation}{0}\renewcommand{\theequation}{9.\arabic{equation}}
\section{Diffusion with Periodical Initial Conditions}
Diffusion process with periodical initial conditions origines for
example in fuel mixing at straight flow engine \cite{Besp}. To
model this process take natural
 $n$ and define markov process
$$Z_n(t)=g(f(y(t)-1+z_n)),\;t\geq 0$$ where r.v. $z_n$ has
uniform distribution on finite set
$$I_n=\left(s={2k+1\over n}:\;\;k=0,1,\ldots,n-1\right)$$
and $z_n,\;y(t)$ are independent. Then random process
$Z_n(t),\;t\geq 0$ may be considered as anomalous diffusion on
interval $[-1,\;1]$ but with periodical initial conditions
\begin{equation}\label{7++}
P(Z_n(0)=-1+s)={1 \over n},\;\;s\in I_n.
\end{equation}

Denote by $\Pi_{t,\;n}=\Pi_{t,\;n}(u),\;\;P_{t,\;n}=P_{t,\;n}(u)$
distribution densities of r.v.`s $Z_n(t),\;f(y(t)-1+z_n),\;\;t>0$.
\begin{lemma}\label{lem3}
The following formula is true
\begin{equation}\label{12}
\Pi_{t,\;n}(u)={1\over n}\sum_{k=0}^{n-1}
\pi_{t,\;1/n}\left(u+1-{2k+1\over n}\right).
\end{equation}
\end{lemma}

\proof Analogously to the formulas (\ref{7}) -- (\ref{5.1+})
obtain
$$ P_{t,\;n}(u)={1 \over n}\sum_{s\in I_n}P_{t}(u+1-s) = {1 \over
n}\sum_{s\in I_n}\sum_{k=-\infty}^{\infty}p_t(u-4k+1-s)
,\;\;-2\leq u< 2,
$$
$$
\Pi_{t,\;n}(u)=\sum_{v:\;g(v)=u}P_{t,\;n}(v)= {1 \over
n}\sum_{s\in
I_n}(\overline{P}_{t}(u+1-s)+\overline{P}_{t}(u+1+s)),\;\;-1\leq
u\leq 1$$ then
\begin{equation}\label{8}
\Pi_{t,\;n}(u)={1 \over
n}\sum_{k=-\infty}^{\infty}p_t\left(u-{2k\over
n}\right),\;\;-1\leq u\leq 1.
\end{equation}
The formula (\ref{8}) leads to
\begin{equation}\label{9}
\Pi_{t,\;n}(u)=\Pi_{t,\;n}\left(u+{2 \over n}\right),\;\;-1\leq
u\leq 1-{2\over n}.
\end{equation}
%В соответствии с формулой (\ref{7+}), определением с.в. $z_n$,
%функции $g$ и множеств $I_n,\;J_n$ имеем
%\begin{equation}\label{10}
%\Pi_{t,\;n}(u)=2P_{t,\;n}(u),\;\;-1\leq u\leq 1.
%\end{equation}
Calculate now the function
$$\Pi_{t,\;n}(u),\;-1\leq u< -1+{2\over n}.$$
For this aim take
$$u=w-1+{1\over n},\;-{1\over n}\leq w< {1\over n}.$$
Then in an accordance with (\ref{5.1}), (\ref{8})
\begin{equation}\label{11}
\Pi_{t,\;n}(w)={\pi_{t,\;1/n}(w)\over n},\;-{1\over n}\leq w<
{1\over n}.
\end{equation}
The formulas (\ref{9}), (\ref{11}) lead to the equality
(\ref{12}). \endproof The equality (\ref{12}) means that the
diffusion (normal or anomalous) on the set $[-1,\;1]$ with
periodical initial conditions (\ref{7++}) and reflecting edges
leads to the same result as a diffusion on isolated (by reflecting
edges) subsets
$$\left[-1+{2k+1 \over n}-{1 \over n},\;-1+{2k+3 \over n}+{1 \over n}\right],\;\;k=0,\ldots,n-1$$
of the set $[-1,\;1].$

Remark that the equality is true for each process $y(t)$ with
independent and symmetrically distributed increments.
\begin{theorem}\label{thm1}
For $tn^{a}\geq 1,\;\;L=[tn^{a}]$:
\begin{equation}\label{14}
||\Pi_{t,\;n}-\pi||\leq 2(1-4Q_1)^{L-1}||P_{1}-P||\;.
\end{equation}
\end{theorem}
\proof The statement of the theorem \ref{thm1} is obtained
directly from the equality (\ref{12}) and the formula (\ref{6.1})
in which $r$ equals to $1/n$.
\endproof
The inequality (\ref{14}) may be interpreted as a decreasing of
characteristic time into $n^{a}$ times in model of anomalous
diffusion with $n-$periodical initial conditions (\ref{7++}):
$$T_n(a,1,1)=T_1(a,1,1)/n^{a}.$$
This result may be easily spread onto general case when r.v.
$Z_n(0)$ has distribution density with continuous derivative
$r_n(u)$ which satisfies periodicity conditions:
$$r_n(u)=r_n\left(u+{2 \over n}\right),\;\;-1\leq u\leq
1-{2\over n}$$ and symmetry conditions
$$r_n(-1+{1\over n}-v)=r_n\left(-1+{1 \over n}+v\right),\;\;0\leq v\leq {1 \over n}$$
and boundary condition
$${d r_n(v)\over dv}=0,\;\;v=\pm 1.$$
It is a generalization of \cite{Besp} results from normal onto
anomalous diffusion.
\setcounter{equation}{0}\renewcommand{\theequation}{10.\arabic{equation}}
\section{Multi-Dimensional Diffusion with Periodic Initial Conditions}
{\bf Multi-Dimensional Anomalous Diffusion.} Results of last
paragraph may be generalized onto two-dimensional case if the
interval $[-1,\;1]$ is replaced by square $[-1,\;1]\times[-1,\;1]$
and process of anomalous diffusion on square with reflecting
boundaries is considered as two independent processes of one
dimension anomalous diffusion on the intervals $[-1,\;1]$ with
reflecting edges. Disintegrate the unit square by rectangular
network with $n^2$ equal squares. Suppose that initial state of
two-dimension diffusion process on unit square has uniform
distribution on the set of these $n^2$ squares centers. Then all
one-dimension estimates of convergence rate are spread onto
two-dimension case.

Denote all one-dimension distributions from last paragraph by
indexes characterized numbers of appropriate coordinates $j=1,2$:
$$P^{(j)}=P^{(j)}(u_j),\;\;P_1^{(j)}=P_1^{(j)}(u_j),\;\;\pi^{(j)}=\pi^{(j)}(u_j),\;\;
\Pi^{(j)}_{t,\;n}=\Pi^{(j)}_{t,\;n}(u_j),\;\;j=1,2.$$ In
accordance with the model of two-dimension diffusion (as a pair of
independent one-dimension diffusion models) obtain:
$$\pi=\pi(u_1,\;u_2)=\pi^{(1)}(u_1)\pi^{(2)}(u_2),\;\;\Pi_{t,\;n}=\Pi_{t,\;n}(u_1,\;u_2)=\Pi^{(1)}_{t,\;n}(u_1)\Pi^{(2)}_{t,\;n}(u_2).$$
From the inequality (\ref{14}) for $L=[tn^{1/a}],\;\;tn^{1/a}\geq
1,\;\;j=1,2$
\begin{equation}\label{14+}
||\Pi^{(j)}_{t,\;n}-\pi^{(j)}||\leq
2(1-4Q_1)^{L-1}||P_{1}^{(j)}-P^{(j)}||=\Delta^{(j)}_{t,\;n}\;,\;\;.
\end{equation}
Define the norm $||\Phi||=\sup\{|\Phi(u_1,\;u_2)|,\;-\infty<
u_1,\;u_2<\infty\}$ of the function $\Phi$ on a plane and put
$\Delta_{t,\;n}=\Delta_{t,\;n}^{(1)}=\Delta_{t,\;n}^{(2)}.$
\begin{theorem}\label{thm2}
For $L=[tn^{1/a}],\;\;tn^{1/a}\geq 1$
\begin{equation}\label{15}
||\Pi_{t,\;n}-\pi||\leq \Delta_{t,\;n}(1+\Delta_{t,\;n}).
\end{equation}
\end{theorem}
\proof From triangle inequality obtain
$$||\Pi_{t,\;n}-\pi||=||\Pi^{(1)}_{t,\;n}\Pi^{(2)}_{t,\;n}-\pi^{(1)}\pi^{(2)}||\leq
||\Pi^{(1)}_{t,\;n}\Pi^{(2)}_{t,\;n}-\pi^{(1)}\Pi^{(2)}_{t,\;n}||+||\pi^{(1)}\Pi^{(2)}_{t,\;n}-\pi^{(1)}\pi^{(2)}||\leq
$$
$$\leq||\Pi^{(2)}_{t,\;n}||\;||\Pi^{(1)}_{t,\;n}-\pi^{(1)}||+||\pi^{(1)}||\;||\Pi^{(2)}_{t,\;n}-\pi^{(2)}||\leq
(||\pi^{(2)}||+||\pi^{(2)}-\Pi^{(2)}_{t,\;n}||)\;||\Pi^{(1)}_{t,\;n}-\pi^{(1)}||+$$
$$+||\pi^{(1)}||\;||\Pi^{(2)}_{t,\;n}-\pi^{(2)}||\leq\left({1\over 2}+||\pi^{(2)}-\Pi^{(2)}_{t,\;n}||\right)
||\Pi^{(1)}_{t,\;n}-\pi^{(1)}||+{1\over
2}\;||\Pi^{(2)}_{t,\;n}-\pi^{(2)}||$$ Using the formula
(\ref{14+}) obtain the inequality (\ref{15}). \endproof As
$\Delta_{t,\;n}\to 0,\;\;t\to\infty$ so convergence rate of
two-dimension case (\ref{15}) is analogous to one-dimension case
(\ref{14}). This result may be easily spread onto multi-dimension
case. Specifics of this multi-dimension diffusion model is that as
one-dimension diffusion components so their initial conditions are
independent.

{\bf Multi-Demensional Normal Diffusion.} Consider $k$-dimension
normal diffusion model in which $\Pi_{t,\;n}(y_1,\ldots,y_k)$ is
distribution density in the moment $t$ at the point
$(y_1,\ldots,y_k)\in[-1,1]^k$:
$$\left(\frac{\partial}{\partial t}-\sum_{j=1}^k\frac{\partial^2}{\partial
y_j^2}\right)\Pi_{t,\;n}(y_1,\ldots,y_k)=0,\;\;$$
$$\frac{\partial}{\partial
y_j}\Pi_{t,\;n}(\pm 1,\ldots,\pm 1,y_j,\pm 1,\ldots,\pm
1)=0,\;\;-1\le y_j\le 1,\;\;j=1,\ldots,k,
$$
$$
\Pi_{0,\;n}(y_1,\ldots,y_k)=\frac{1}{2^k}+\sum_{j_1,\ldots,j_k=1}^{m}a(j_1,\ldots,j_k)\prod_{r=1}^k
\cos \pi nj_r y_{r}>0,\;\;m<\infty.$$ Then
$$
\Pi_{t,\;n}(y_1,\ldots,y_k)=\frac{1}{2^k}+\sum_{j_1,\ldots,j_k=1}^{m}a(j_1,\ldots,j_k)\prod_{r=1}^k
\cos \pi nj_r y_{r}\exp\left(-\pi^2n^2 t\sum_{j=1}^m
k_j^2\right)$$ и при $|a(1,\ldots,1)|>0$
$$
||\Pi_{t,\;n}(y_1,\ldots,y_k)-\frac{1}{2^k}||\sim|a(1,\ldots,1)|\exp(-\pi^2
k n^2 t ),\;t\to \infty.
$$
So it is possible to choose $n$-periodic (by each coordinate)
initial conditions that
$$T_n(2,1,k)=T_1(2,1,1)/kn^{2}.$$

The paper is supported by RFBR, projects 03-01-00512,
03-07-07-90334v.

\end{document}